----------
X-Sun-Data-Type: default
X-Sun-Data-Description: default
X-Sun-Data-Name: tadek4
X-Sun-Charset: us-ascii
X-Sun-Content-Lines: 405

\input amstex \documentstyle{amsppt}
\magnification 1200 \def\R{\text{\bf R}} \def\sF{\Cal F}
\hoffset= 0.0in
\voffset= 0.0in
\hsize=32pc
\vsize=43pc
\baselineskip=24pt
\NoBlackBoxes
\topmatter
\author A. Dranishnikov and T. Januszkiewicz \endauthor
\title Every Coxeter group acts amenably on a compact space \endtitle
\date{July 10, 1999} \enddate
\abstract 
Coxeter groups are Higson-Roe amenable, i.e. they admit amenable actions
on compact spaces.
Moreover, they have finite asymptotic dimension.
This answers affirmatively a question from [H-R].
\endabstract
\thanks 
The first author was partially supported by NSF grant DMS-9971709.
The second author was partially supported by KBN grant.
2 P03A 023 14. Both authors enjoyed hospitality of IHES
while this paper was conceived.
\endthanks
\address Pennsylvania State University, Department of Mathematics, 
218 McAllister Buil\-ding, University Park, PA 16802 , USA
\endaddress
\email  dranish\@math.psu.edu
\endemail
\address Instytut Matematyczny PAN \hbox{ and } Instytut Matematyczny 
Uniwersytetu Wroc\-\l aw\-skie\-go, pl. Grunwaldzki 2/4; 50-384 Wroc\l aw,
Poland. \endaddress
\email  tjan\@math.uni.wroc.pl \endemail
\endtopmatter
\document
{\bf 1. Higson-Roe amenability.}
An action of a discrete group $G$ on a compact space $X$ is {\it topologically amenable}
[AD-R] if there is a sequence of continuous maps $b^n:X\to P(G)$
to the space of probability measures on $G$ with the weak$^*$-topology such that for every $g\in G$,
$\lim_{n\to\infty}\sup_{x\in X}\|gb^n_x-b_{gx}^n\|_1=0$.
Here a measure $b^n_x=b^n(x)$ is considered as a function $b^n_x:G\to[0,1]$ and
$\|\ \|_1$ is the $l_1$-norm.
\proclaim{Definition} 
A discrete countable group $G$ is called  Higson-Roe amenable if $G$ admits
a topologically amenable action on a compact space or, equivalently, 
its natural action
on the Stone-\v{C}ech compactification $\beta G$ is topologically amenable.
\endproclaim
\bigskip{\bf 2. Property A. } 
The property A was introduced in [Yu]. 
For metric spaces of bounded geometry
it was reformulated in [H-R] as follows.
\proclaim{Property A} A discrete metric space $Z$ has the property $A$ if and
only if there is a sequence of maps $a^n:Z\to P(Z)$ such that
\roster
\item{} for every $n$ there is some $R>0$ with the property that for
every $z\in Z,\ \ \ 
supp(a^n_z)\subset\{z'\in Z\mid d(z,z')<R\}$
and
\item{} for every $K>0$, $\lim_{n\to\infty}\sup_{d(z,w)<K}\|a^n_z-a^n_w\|_1=0$.
\endroster
\endproclaim
\proclaim{Lemma [H-R]}
A finitely generated group $G$ is Higson-Roe amenable if and only if 
the underlying
metric space $G$ with a word metric has property $A$.
\endproclaim

A tree $T$ posesses a natural metric where 
every edge has the length one. We denote by $V(T)$ the set of vertices of $T$
with induced metric. The idea of the proof of the following proposition is
taken from [Yu].
\proclaim{Proposition 1}
For any tree $T$ the metric space $V(T)$ has property $A$.
\endproclaim
\demo{Proof}
Let $\gamma_0:\R\to T$ be a geodesic ray in $T$, i.e. an isometric 
embedding of the
half-line $\R$. For every point $z\in V(T)$ there is a unique geodesic
ray $\gamma_z$ issued from $z$ which intersects $\gamma_0$ 
along a geodesic ray. Let $V=V(T)\cap[z,\gamma_z(n)]\subset im(\gamma_z)$.
 We define $a^n_z$ as a Dirac measure
$\Sigma_{v\in V}\frac{1}{n+1}\delta_v$
supported uniformly by vertices lying in the geodesic segment $[z,\gamma_z(n)]$.
Then condition (1) holds. If $d(z,w)<K$, then the geodesic segments
$[z,\gamma_z(n)]$ and $[w,\gamma_w(n)]$ overlap along a geodesic segment of 
the length $\ge n-2K$. Then $\sup\|a^n_z-a^n_w\|_1\le \frac{2K}{n+1}$ and hence
the condition (2) holds.\qed
\enddemo

{\bf Remark. } 
We do not assume that a tree is 
locally finite in the above proposition.
Thus the geometry is not bounded, and 
the variant of the property A we use may differ from the
genuine property A of [Yu].

The following proposition is obvious.
\proclaim{Proposition 2}
Let $Z=Z_1\times Z_2$ be the product of discrete metric spaces with the
$l_1$ metric. Assume that $Z_1$ and $Z_2$ have property A. 
Then $Z$ has property A.
\endproclaim
\proclaim{Proposition 3}
Let $W\subset Z$ and let $Z$ have property A. Then $W$ has property $A$.
\endproclaim
\demo{Proof}
Let $a^n:Z\to P(Z)$ be a sequence of maps from the definition of property
A for the space $Z$. We define a sequence $A^n:W\to P(W)$ by the following 
rule.
Let $r:Z\to W$ be a retraction which takes every point $z\in Z$ to a closest
point $r(z)\in W$. We define $A^n_w=P(r)(a^n_w)$ where $P(r):P(Z)\to P(W)$ is
the induced map on probability measures. If the support $supp(a^n_w)$ lies
in the ball $B_R(w)\subset Z$, then $supp(A^n_w)$ lies in the ball $B_{2R}(w)$.
Therefore the condition (1) holds. 
Let $a^n_z=\Sigma_{x\in Z}\lambda_x\delta_x$ and 
$a^n_w=\Sigma_{x\in Z}\nu_x\delta_x$. 
Note that $\|A^n_z-A^n_w\|_1=\|P(r)(a^n_z)-P(r)(a^n_w)\|_1=
\|\Sigma_{x\in Z}\lambda_x\delta_{r(x)}-\Sigma_{x\in Z}\nu_x\delta_{r(x)}\|_1=
\Sigma_{y\in W}|\Sigma_{z\in r^{-1}(y)} (\lambda_z-\nu_z)|\le
\Sigma_{x\in Z}|\lambda_z-\nu_z|=\|a^n_z-a^n_w\|_1$. Hence the condition 
(2) holds.
\qed
\enddemo

\bigskip
{\bf 3. Coxeter groups.}
Recall that a Coxeter group 
$(\Gamma,W)$ is a group $\Gamma$ 
with a distinguished set of generators 
$w_i\in W$  and relations $w_i^2=1=(w_iw_j)^{m_{ij}}$, where $m_{ij}$ is 
zero (and then there is no relation between $w_i$ and $w_j$) or an integer 
$\geq2$.

\medskip
For any Coxeter group $(\Gamma,W)$ there is a cell complex $C(\Gamma)$ on 
which $\Gamma $ acts properly. Its construction, due to M. Davis, %
%
%
%
%
is described fully in [D].
It is defined as follows: cells are indexed by right cosets $\Gamma/\Gamma_S$,
where $\Gamma_S$  is  a {\it finite} group generated by a subset $S$ of $W$. 
A cell $[\gamma]$ is a face of $[\eta]$ if $[\gamma]\subset[\eta]$ as cosets. 
Its vertices correspond to elements of $C(\Gamma)$, edges are indexed by 
generators etc. Its 1-skeleton is the Cayley graph of $(\Gamma, W)$.

The obvious action of $\Gamma$ on  $C(\Gamma)$ coming from the left action
of $\Gamma $ on itself is a reflection group action. Any reflection 
(i.e. an element conjugated to a generator in $W$) has its mirror of fixpoints.
Any mirror is two sided i.e. its complement has two components. 
Closures of connected components of the set of all mirrors are called 
fundamental domains.

\proclaim{Theorem 1 [J]}
Every Coxeter group $\Gamma$ can be isometrically imbedded
in a finite product of trees $\Pi T_i$ with $l_1$ metric on it in such a way that
the image of $\Gamma$ under this embedding is contained in the set of
vertices of $\Pi T_i$.
\endproclaim

\proclaim{Lemma 1 (\cite M)}
Let $\Gamma_0$ be a normal torsion free subgroup of $\Gamma$.
Then for each mirror $H$ and every $\gamma_0\in \Gamma_0$,
$H\cap \gamma_0(H)$ is either $H$ or the empty set.
\endproclaim

Observe that such $\Gamma_0$, even with an additional property
of being of finite index in $\Gamma$, exists by the Selberg Lemma.
The proof of the lemma is just two lines, so we repeat it. 

\demo{Proof}
Let $h$ be the reflection in $H$, and consider the product of reflections
$g=h\gamma_0h\gamma_0^{-1}$. If $H\cap \gamma_0(H)$ is nonempty,
$g$ is a torsion element. On the other hand since $\Gamma_0$ is normal, 
$g$ is in $\Gamma_0$, hence it is an identity.
Thus $h$ commutes with $\gamma_0$, and $H = \gamma_0(H)$.
$\square$
\enddemo

Let $\Cal H$ be the set of orbits for  $\Gamma_0$ action on the set of all 
mirrors. Fix an $h\in \Cal H$. Define a graph $T_h$ as follows. 
Its vertices are connected components of $C(\Gamma)-\bigcup_s \Gamma_0 (s)$ 
where $s$ is any mirror in $h$.
Two vertices are joined by an edge if the components are adjacent in 
$C(\Gamma)$ i.e. if they intersect after taking closures. 

\proclaim {Lemma 2} $T_h$ is a tree. \endproclaim
\demo{Proof} Any loop $\lambda$ in $T_h$ lifts to a path in $C(\Gamma)$
which can be closed up to a loop $\Lambda$ without crossing the $h$--mirror.
Projection of $\Lambda$ is again $\lambda$.
Since $C(\Gamma)$ is contractible (cf. \cite D), $\Lambda$, hence $\lambda$,
are homologous to zero thus $T_h$ is a tree. $\square$
\enddemo

There is an obvious simplicial map from $C(\Gamma)$ to $T_h$,  given on 
vertices by $g \rightarrow [g]_h$, 
that is mapping a fundamental domain to the 
connected component of $C(\Gamma)-\bigcup_s \Gamma_0 (s)$ it belongs to.
Clearly this map is  $\Gamma_0$ equivariant. We take the diagonal of the 
family $\mu :C(\Gamma)\to \prod_h T_h$ to get a $\Gamma_0$ 
equivariant embedding of the Davis complex into the product of trees.

\proclaim{Lemma 3}
The map $\mu$ is a $\Gamma$--equivariant embedding
\endproclaim
\demo{Proof}
First notice that $\Gamma$, in fact $\Gamma\over\Gamma_0$,
acts on $\Cal H$. An element $g$ 
maps the tree $T_h$ to $T_{g(h)}$ simplicially.
Thus $\Gamma$ acts on $\prod_h T_h$  by permuting factors of the product.
Explicitly
$g(x_{h_1} ,\dots, x_{h_n})=
(gx_{g^{-1}(h_1)} ,\dots, gx_{g^{-1}(h_n)})$
Now equivariance of $\mu$ is obvious.
To see that it is an embedding, notice that two  vertices of $C(\Gamma)$
differ iff they are separated by some mirror say in $h$; thus their images 
in $T_h$ are different.$\square$
\enddemo

\demo {Proof of Theorem 1}
In view of Lemmas 1--3 we should only notice that teh $l_1$--metric
on the product of trees restricted to
the image of $\mu$ agrees with the word metric on $\Gamma$.
This is clear, as the word metric $d(g,h)$ counts the number
of mirrors between chambers $gF$ and $hF$ for any fixed chamber $F$
in the Davis complex.
\enddemo

\proclaim{Theorem A}
Every Coxeter group $\Gamma$ is Higson-Roe amenable.
\endproclaim
\demo{Proof}
In view of Higson-Roe Lemma it suffices to show that $\Gamma$ has property A.
Propositions 1, 2, 3 and Theorem 1 imply the result.
\qed\enddemo

\bigskip
{\bf 4. Asymptotic dimension.}
The following definition is analogous to Ostrand's characterization of
covering dimension.
\proclaim{Definition [Gr]}
The asymptotic dimension of a metric space $X$ does not acceed $n$,
$as\dim X\le n$ if for arbitrary large $d>0$ there are $n+1$ uniformly bounded
$d$-disjoint families $\sF_i$ of sets in $X$ such that the union $\cup\sF_i$
forms a cover of $X$. A family $\sF$ is $d$-disjoint provided 
$\min\{dist(x,y)\mid x\in F_1,\
y\in F_2,\ F_1\ne F_2,\ \ F_1,F_2\in\sF\}\ge d$.
\endproclaim
\proclaim{Proposition 4}
For every tree $T$, \  $as\dim T\le 1$.
\endproclaim
\demo{Proof}
Let $x_0\in T$ be a fixed point and let $B_{nd}(x_0)$ denote the closed ball
of radius $nd$ centered at $x_0$. Let $R_{n,i}=
B_{nd}\setminus Int(B_{(n-i)d})$ .
We define 

$\sF_1=\{\text{components of}\ R_{n,2}\ \ \text{intersected with}
\ R_{n,1}\ \ \text{for} \ n \ \text{odd}\}$
and 

$\sF_2=\{\text{components of}\ R_{n,2}\ 
\ \text{intersected with}\ \ R_{n,1}
 \ \text{for} \ n \ \text{even}\}$.
The diameter of each component of $R_{n,2}$ is less than $4d$ and for every two components
in $R_{n,2}$ the distance between their traces on $R_{n,1}$ is greater than $2d$.\qed
\enddemo
The following proposition follows from the reformulation of the definition
of asymptotic dimension in terms of anti-\v Cech approximation of metric
spaces by polyhedra [Gr].
\proclaim{Proposition 5}
$as\dim(X\times Y)\le as\dim X + as\dim Y$.
\endproclaim
\proclaim{Theorem B}
For every Coxeter group $\Gamma$, $as\dim\Gamma<\infty$.
\endproclaim
\demo{Proof}
Theorem 1 and Propositions 4, 5 imply the proof.
\qed
\enddemo

{\bf Remarks. }

1. Theorem B and Lemma 4.3 of [H-R] give an alternative proof of
Theorem A.

2. Any of the Theorems A and B implies that Coxeter groups
admit a coarsly uniform embedding into Hilbert space. This
result is due to [B-J-S].

3. We would like to point out (but not go into details)
that the proof we presented for both theorem A and B works
in much greater generality: that of "zonotopal complexes of 
nonpositive curvature" ([D-J-S]), perhaps with an additional property of 
"foldability". An interesting subclass here is that of cubical 
nonpositively curved complexes.

\enddocument